\documentclass[12pt,reqno,a4wide]{amsart}

\usepackage{tikz} 
\usepackage{calrsfs}
\usepackage{mathrsfs}

\usepackage{array}
\usepackage{hyperref}

\usetikzlibrary{decorations.pathreplacing}
\usetikzlibrary{fit}

\usepackage{pgfplots}

\allowdisplaybreaks

\begin{document}
\bibliographystyle{plain}

%
%

	\title [Partial Dyck path interpretation]
	{Partial Dyck path interpretation for three sequences   in the Encyclopedia of Integer Sequences}

	\author[H. Prodinger ]{Helmut Prodinger }
	\address{Department of Mathematics, University of Stellenbosch 7602, Stellenbosch, South Africa
	and
NITheCS (National Institute for
Theoretical and Computational Sciences), South Africa.}
	\email{hproding@sun.ac.za}

	
	\begin{abstract}
Descents of odd length in Dyck paths are discussed, taking care of some variations. 
The approach is based on generating functions and the kernel method 
 and augments relations about them
from the Encyclopedia of Integer Sequences, that were pointed out by David Callan \cite{Callan}.
	\end{abstract}
	

\maketitle

\section{Introduction}
A Dyck path consists of (a sequence of) up-steps and down-steps of one unit each, such that it never goes below the $x$-axis and ends at the $x$-axis.
There are numerous entries of this popular concept in \cite{OEIS}. If the path does not end on the $x$-axis but can be completed to a Dyck path, we
speak about a partial Dyck path or the prefix of a Dyck path.  A \emph{descent} is a maximal sequence of contiguous down-steps; its length is measured in
terms of down-steps. We follow Callan \cite{Callan} and interpret the sequences A101785, A113337, and A143017 in \cite{OEIS}; however we extend the analysis by considering partial versions of such paths. As a bonus, we discuss a further sequence of related interest. 

The subclasses of Dyck paths are:
\begin{itemize}
	\item All descents have odd length,
	\item All descents have odd length, but the last descent has even length,
		\item All descents have odd length, but the last descent can have even or odd length.
\end{itemize}
Since we do not consider the empty path (of length zero), it always makes sense to speak about the ``last descent.''

Our bonus is about descents of odd length if they do not end on the $x$-axis, but any descent reaching the $x$-axis can be even resp.\ odd (not only the final one).

\section{All descents have odd length}
We use an (infinite) graph (automaton) to control the odd length(s) of the descents:
  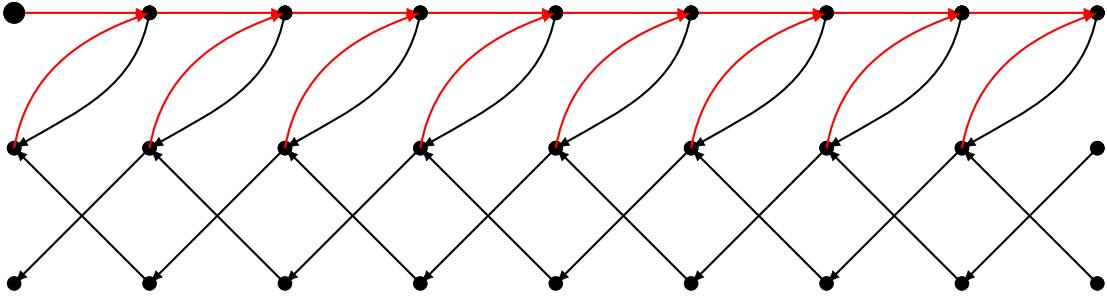
\begin{figure}[h]
  	\label{F1}

 	\begin{center}
 		\begin{tikzpicture}[scale=1.8,main node/.style={circle,draw,font=\Large\bfseries}]

 			\foreach \x in {0,1,2,3,4,5,6,7,8}
 			{
 				\draw (\x,0) circle (0.05cm);
 				\fill (\x,0) circle (0.05cm);
 				\draw (\x,-1) circle (0.05cm);
 				\fill (\x,-1) circle (0.05cm);
 			}

 			\fill (0,0) circle (0.08cm);

 			\foreach \x in {0,1,2,3,4,5,6,7}
 			{
 				\draw[ thick,red, latex-] (\x+1,0) to  (\x,0);	
\draw[ thick,  latex-] (\x,-2)   to  (\x+1,-1);
\draw[ thick,  latex-] (\x,-1) to  (\x+1,-2);
 			}

 		\foreach \x in {1,2,3,4,5,6,7}
 	{
 	}



 			\foreach \x in {0,1,2,3,4,5,6,7,8}
{
	\draw (\x,0) circle (0.05cm);
	\fill (\x,0) circle (0.05cm);
	\draw (\x,-1) circle (0.05cm);
	\fill (\x,-1) circle (0.05cm);
}

\foreach \x in {1,2,3,4,5,6,7,8}
{
	\draw (\x,0) circle (0.05cm);
	\fill (\x,0) circle (0.05cm);
}

	\foreach \x in {0,1,2,3,4,5,6,7,8}
{
	\draw (\x,0) circle (0.05cm);
	\fill (\x,0) circle (0.05cm);
	\draw (\x,-2) circle (0.05cm);
	\fill (\x,-2) circle (0.05cm);
}

\fill (0,0) circle (0.08cm);

\foreach \x in {0,1,2,3,4,5,6,7}
{ 				
	\draw  [thick, red, -latex] (\x,-1)[in=200, out=80] to  (\x+1,0);	
	\draw  [thick,latex-] (\x,-1)[in=-100, out=30] to  (\x+1,0);	
}

 		\end{tikzpicture}
 	\caption{Three layers of states, labelled $f,g,h$, in that order. The state $g_0$ is responsible for allowed Dyck paths.}
 	\end{center}
 \end{figure}
The recursions for the generating functions $f_i=f_i(z)$, $g_i=g_i(z)$, $h_i=h_i(z)$,  can be read off from the graph in Fig.~\ref{F1}, by considering the last step made;
\begin{align*}
f_{i+1}&=zf_{i}+zg_{i},\ i\ge 0,\ f_0=1, \\
g_{i}&=zf_{i+1}+zh_{i+1},\ i\ge0,\\
h_{i}&=zg_{i+1},\ i\ge0.
\end{align*}
So the coefficient $[z^n]f_i$ is the number of partial Dyck path under consideration that end in the $i$-th state of the first layer, and similarly for the other two.
Summing up the recursions we get
\begin{align*}
	\sum_{i\ge0}u^{i+1}f_{i+1}&=\sum_{i\ge0}u^{i+1}zf_{i}+\sum_{i\ge0}u^{i+1}zg_{i}, \\
	\sum_{i\ge0}u^{i+1}g_{i}&=\sum_{i\ge0}u^{i+1}zf_{i+1}+\sum_{i\ge0}u^{i+1}zh_{i+1},\\
	\sum_{i\ge0}u^{i+1}h_{i}&=\sum_{i\ge0}u^{i+1}zg_{i+1}
\end{align*}
and with
\begin{equation*}
F(u)=\sum_{i\ge0}u^if_i,\quad G(u)=\sum_{i\ge0}u^ig_i,\quad H(u)=\sum_{i\ge0}u^ih_i,
\end{equation*}
furthermore
\begin{align*}
	F(u)&=1+uzF(u)+uzG(u), \\
	uG(u)&=zF(u)-z+zH(z)-zh_0,\\
	uH(u)&=zG(u)-zg_0.
\end{align*}
This system can be solved:
\begin{align*}
F(u)&=1+\frac{uz(uzh_0+z^2g_0-u^2+z^2)}{u^3z+u^2z^2-u^2-z^3u+z^2},\\
G(u)&=\frac{u(-u^2z-u^2zh_0+uh_0-uz^2g_0+zg_0)}{u^3z+u^2z^2-u^2-z^3u+z^2},\\
H(u)&=\frac{u(-g_0zu^2+g_0u-z^2h_0u-z^2g_0u-z^2u+zh_0)}{u^3z+u^2z^2-u^2-z^3u+z^2}.
\end{align*}
Plugging in $u=0$ does not help, but the (joint) denominator may be factored:
\begin{equation*}
u^3z+u^2z^2-u^2-z^3u+z^2=z(u-v_1)(u-v_2)(u-v_3).
\end{equation*}
Two of these factors are `bad,' as the have no combinatorial relevance and, more importantly, no power series expansion around $(u,z)=(0,0)$. 
This procedure is commonly known as the \emph{kernel method}, see e.g. \cite{prodinger-kernel}. Thus, 
$(u-v_2)(u-v_3)$ may be divided out from both, numerator and denominator, yielding
\begin{align*}
F(u)&=1-\frac{u}{u-v_1},\\
	G(u)&=\frac{-z(1+h_0)}{u-v_1},\\
	H(u)&=\frac{-zg_0}{u-v_1}
\end{align*}
and thus
\begin{equation*}
g_0=\frac{zv_1}{-z^2+v_1^2},\quad h_0=\frac{z^2}{-z^2+v_1^2}=-\frac{v_1}{z}-1+\frac1{z^2}.
\end{equation*}
The relevant solution can be expanded by a Computer algebra system like Maple:
\begin{equation*}
v_1=\frac1z-z-z^5-2z^7-4z^9-10z^{11}-26z^{13}-68z^{15}-183z^{17}-504z^{19}-1408z^{21}-\dots
\end{equation*}
As predicted in OEIS (and automatically established by Computer algebra)
\begin{equation*}
g_0=\frac{z^2(1+g_0)}{1-z^4(1+g_0)^2}
\end{equation*}
and thus we get the algebraic equation (with $z^2=Z$)
\begin{equation*}
-{Z}^{2} g_0^{3}-2{Z}^{2} g_0^{2}+{ g_0}-{Z}^{2}{ g_0}-Z{ g_0}-Z=0.
\end{equation*}
Note that 
\begin{equation*}
g_0=Z+Z^2+2Z^3+5Z^4+12Z^5+30Z^6+79Z^7+213Z^8+584Z^9+1628Z^{10}+4600Z^{11}+\dots,
\end{equation*}
which is sequence A101785 in \cite{OEIS}.

Reading off the coefficient of $u^j$ in $F(u)$, $G(u)$, $H(u)$ gives us explicit expressions for $f_j$, $g_j$, $h_j$ and thus all types of \emph{partial} Dyck paths under consideration. We have $f_0=1$ and for $j\ge1$
\begin{equation*}
[u^j]F(u)=-[u^j]\frac{u}{u-v_1}=[u^{j-1}]\frac{1}{v_1(1-u/v_1)}=v_1^{-j}.
\end{equation*}
Further,
\begin{equation*}
	[u^j]G(u)=[u^j]\frac{-z(1+h_0)}{u-v_1}=[u^j]\frac{z(1+h_0)}{v_1(1-u/v_1)}=\frac{z(1+h_0)}{v_1^{j+1}};
\end{equation*}
finally
\begin{equation*}
	[u^j]H(u)=[u^j] \frac{-zg_0}{u-v_1}=[u^j] \frac{zg_0}{v_1(1-u/v_1)}=\frac{zg_0}{v_1^{j+1}}.
\end{equation*}

\section{All descents have odd length, but the last descent has even length}

For this instance we don't need  new computations, as it is covered by $h_0$ from the previous section:
\begin{equation*}
h_0=Z^2+2Z^3+4Z^4+10Z^5+26Z^6+68Z^7+183Z^8+504Z^9+1408Z^{10}+3982Z^{11}+\dots
\end{equation*}
This is sequence  A113337.

\section{All descents have odd length, but the last descent can have even or odd length}

Again, no new computation is necessary, as we just have to consider $g_0+h_0$;
\begin{equation*}
Z+2 Z^2+4 Z^3+9 Z^4+22 Z^5+56 Z^6+147 Z^7+396 Z^8+1088 Z^9+3036 Z^{10}+8582 Z^{11}+\dots,
\end{equation*}
which is sequence A143017.

\section{Bonus problem: descents of odd length when they do not touch the $x$-axis, otherwise no restriction}

In this instance, we draw again (Figure~\ref{F2}) a graph and notice the anomalies in the beginning (related to returns to the $x$-axis).
We use the same letters as before, but now with a slightly different meaning. As a first step, it is beneficial to ignore   $f_0$, $g_0$, $h_0$ 
for the time being, which leads us to Figure~\ref{F3}.

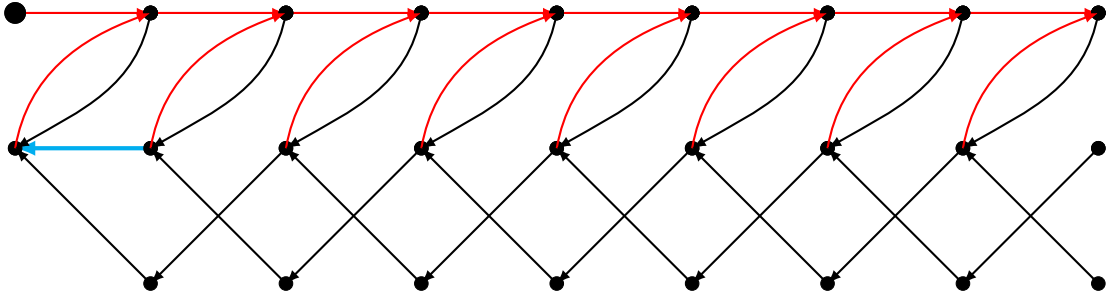
\begin{figure}[h]
\label{F2}	
	
	\begin{center}
		\begin{tikzpicture}[scale=1.8,main node/.style={circle,draw,font=\Large\bfseries}]

			\foreach \x in {0,1,2,3,4,5,6,7,8}
			{
				\draw (\x,0) circle (0.05cm);
				\fill (\x,0) circle (0.05cm);
				\draw (\x,-1) circle (0.05cm);
				\fill (\x,-1) circle (0.05cm);
			}

			\fill (0,0) circle (0.08cm);

			\foreach \x in {0,1,2,3,4,5,6,7}
			{
				\draw[ thick,red, latex-] (\x+1,0) to  (\x,0);	
				\draw[ thick,  latex-] (\x,-1) to  (\x+1,-2);
			}			
			
			\foreach \x in {1,2,3,4,5,6,7}
			{
				\draw[ thick,  latex-] (\x,-2)   to  (\x+1,-1);
			}
			
			\foreach \x in {1,2,3,4,5,6,7}
			{
			}			
			
							\draw[ ultra thick, cyan,  latex-] (0,-1)   to  (1,-1);


			
			\foreach \x in {0,1,2,3,4,5,6,7,8}
			{
				\draw (\x,0) circle (0.05cm);
				\fill (\x,0) circle (0.05cm);
				\draw (\x,-1) circle (0.05cm);
				\fill (\x,-1) circle (0.05cm);
			}
			
			\foreach \x in {1,2,3,4,5,6,7,8}
			{
				\draw (\x,0) circle (0.05cm);
				\fill (\x,0) circle (0.05cm);
			}
			
			\foreach \x in {0,1,2,3,4,5,6,7,8}
			{
				\draw (\x,0) circle (0.05cm);
				\fill (\x,0) circle (0.05cm);
			}

			\foreach \x in {1,2,3,4,5,6,7,8}
{
	\draw (\x,0) circle (0.05cm);
	\fill (\x,0) circle (0.05cm);
	\draw (\x,-2) circle (0.05cm);
	\fill (\x,-2) circle (0.05cm);
}

			\fill (0,0) circle (0.08cm);
			
			\foreach \x in {0,1,2,3,4,5,6,7}
			{ 				
				\draw  [thick, red, -latex] (\x,-1)[in=200, out=80] to  (\x+1,0);	
				\draw  [thick,latex-] (\x,-1)[in=-100, out=30] to  (\x+1,0);	
			}

		\end{tikzpicture}
		\caption{Three layers of states, labelled $f,g,h$, in that order. }
	\end{center}
\end{figure}
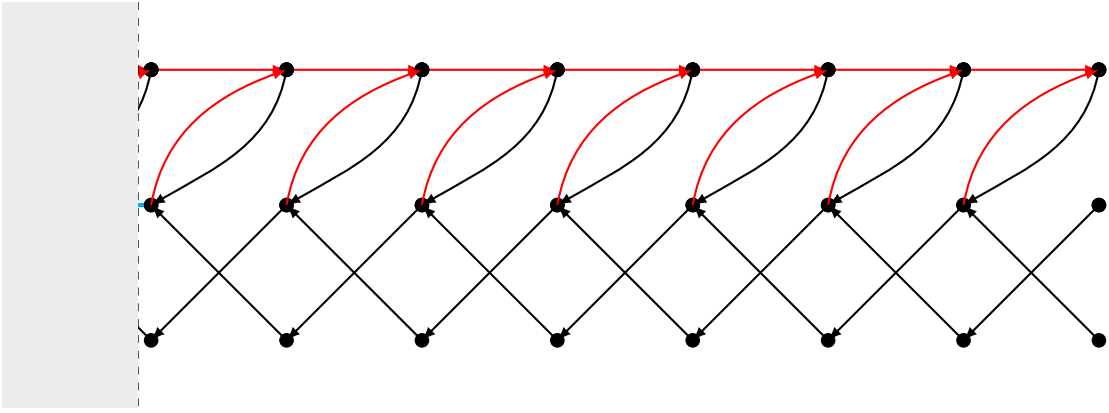
\begin{figure}[h]
	\label{F3}
	
	\begin{center}
		\begin{tikzpicture}[scale=1.8,main node/.style={circle,draw,font=\Large\bfseries}]

			\foreach \x in {0,1,2,3,4,5,6,7,8}
			{
				\draw (\x,0) circle (0.05cm);
				\fill (\x,0) circle (0.05cm);
				\draw (\x,-1) circle (0.05cm);
				\fill (\x,-1) circle (0.05cm);
			}

			\fill (0,0) circle (0.08cm);

			\foreach \x in {0,1,2,3,4,5,6,7}
			{
				\draw[ thick,red, latex-] (\x+1,0) to  (\x,0);	
				\draw[ thick,  latex-] (\x,-1) to  (\x+1,-2);
			}			
			
			\foreach \x in {1,2,3,4,5,6,7}
			{
				\draw[ thick,  latex-] (\x,-2)   to  (\x+1,-1);
			}
			
			\foreach \x in {1,2,3,4,5,6,7}
			{
			}			
			
			\draw[ ultra thick, cyan,  latex-] (0,-1)   to  (1,-1);


			
			\foreach \x in {0,1,2,3,4,5,6,7,8}
			{
				\draw (\x,0) circle (0.05cm);
				\fill (\x,0) circle (0.05cm);
				\draw (\x,-1) circle (0.05cm);
				\fill (\x,-1) circle (0.05cm);
			}
			
			\foreach \x in {1,2,3,4,5,6,7,8}
			{
				\draw (\x,0) circle (0.05cm);
				\fill (\x,0) circle (0.05cm);
			}
			
			\foreach \x in {0,1,2,3,4,5,6,7,8}
			{
				\draw (\x,0) circle (0.05cm);
				\fill (\x,0) circle (0.05cm);
			}
			
			\foreach \x in {1,2,3,4,5,6,7,8}
			{
				\draw (\x,0) circle (0.05cm);
				\fill (\x,0) circle (0.05cm);
				\draw (\x,-2) circle (0.05cm);
				\fill (\x,-2) circle (0.05cm);
			}

			\fill (0,0) circle (0.08cm);
			
			\foreach \x in {0,1,2,3,4,5,6,7}
			{ 				
				\draw  [thick, red, -latex] (\x,-1)[in=200, out=80] to  (\x+1,0);	
				\draw  [thick,latex-] (\x,-1)[in=-100, out=30] to  (\x+1,0);	
			}
			
\draw[thick, dashed] (1-0.1,0.5) to (1-0.1,-2.5);	
\path[fill=lightgray!30, fill opacity=0.6]	(1-0.1,0.5) to (1-0.1,-2.5) to 	(0-0.1,-2.5) to (0-0.1,0.5) ;
			
		\end{tikzpicture}
		\caption{Three layers of states, labelled $f,g,h$, in that order. }
	\end{center}
\end{figure}
It is beneficial to consider first the following generating functions,
\begin{equation*}
F(u)=\sum_{i\ge1}u^{i-1}f_i,\quad G(u)=\sum_{i\ge1}u^{i-1}g_i,\quad H(u)=\sum_{i\ge1}u^{i-1}h_i
\end{equation*}
following the recursions
\begin{align*}
f_{i+1}&=zf_i+zg_i,\ i\ge1,\\
 g_i&=zf_{i+1}+zh_{i+1}	,\ i\ge1,\\
 h_i&=zg_{i+1}	,\ i\ge1.
\end{align*}
Then
\begin{align*}
	\sum_{i\ge1}u^if_{i+1}&=\sum_{i\ge1}u^izf_i+\sum_{i\ge1}u^izg_i,\\
	\sum_{i\ge1}u^ig_i&=\sum_{i\ge1}u^izf_{i+1}+\sum_{i\ge1}u^izh_{i+1},\\
	\sum_{i\ge1}u^ih_i&=\sum_{i\ge1}u^izg_{i+1}
\end{align*}
and
\begin{align*}
	F(u)-f_1&=uzF(u) +uzG(u),\\
	uG(u)&=zF(u)-zf_1+zH(u)-zh_1,\\
	uH(u)&=zG(u)-zg_1.
\end{align*}
The system is almost the same as in previous sections except for the unspecified quantity $f_1$ which will be fixed later;
\begin{align*}
	F(u)&=f_1+\frac{ uz(-f_1u^2+zh_1u+z^2g_1+z^2f_1)}{u^3z+u^2z^2-u^2-z^3u+z^2},\\
	G(u)&=\frac{z(uh_1+zg_1-u^2zf_1-u^2zh_1-uz^2g_1)}{u^3z+u^2z^2-u^2-z^3u+z^2},\\
	H(u)&=\frac{z(zh_1-uz^2f_1-uz^2h_1-uz^2g_1+g_1u-g_1u^2z)}{u^3z+u^2z^2-u^2-z^3u+z^2}.
\end{align*}
The denominator is the same as before, and dividing out $(u-v_2)(u-v_3)$ leads to
\begin{align*}
	F(u)&=f_1-\frac{uf_1}{u-v_1},\\
	G(u)&= \frac{-z(f_1+h_1)}{u-v_1},\\
	H(u)&=\frac{-zg_1}{u-v_1}.
\end{align*}
Consequently, setting $u=0$,
\begin{align*}
	g_1&=\frac{zf_1v_1}{-z^2+v_1^2},\quad
h_1=\frac{z^2f_1}{-z^2+v_1^2}.
\end{align*}
Now, bringing the initial layers into the game, we have
\begin{equation*}
f_1=z+zg_0,\quad g_0=z(f_1+g_1+h_1),
\end{equation*}
and therefore the quantity of interest
\begin{equation*}
g_0=\frac{z^2v_1}{v_1(1-z^2)-z}.
\end{equation*}
The series expansion is
\begin{equation*}
g_0=Z+2 Z^2+5 Z^3+13 Z^4+35 Z^5+97 Z^6+274 Z^7+785 Z^8+2275 Z^9+6655 Z^{10}+19618 Z^{11}+\dots
\end{equation*}
The series $g_0$ satisfies a cubic equation, and all functions $[u^j]F(u)$, $[u^j]G(u)$, $[u^j]H(u)$ could be computed as well.
\bibliographystyle{plain}


\end{document}